\documentclass{article}
\usepackage{amsmath, amsfonts, amsthm, amssymb}
\usepackage{graphicx}
\usepackage{float}
\usepackage{verbatim}
\usepackage{color}
\usepackage[sort]{cite}

\numberwithin{equation}{section}

\newtheorem{thm}{Theorem}[section]
\newtheorem{lma}[thm]{Lemma}
\newtheorem{cor}[thm]{Corollary}

\newtheorem{conj}[thm]{Conjecture}

\renewcommand{\geq}{\geqslant}
\renewcommand{\leq}{\leqslant}

\allowdisplaybreaks

\title{The Assouad spectrum and the quasi-Assouad dimension: \\ a tale of two spectra}

\author{Jonathan M. Fraser$^*$, Kathryn E. Hare$^\dagger$, Kevin G. Hare$^\dagger$, Sascha
Troscheit$^\dagger$, Han Yu$^*$}

\begin{document}

\maketitle

\vspace{-5mm}
\begin{center}
	$^*$Mathematical Institute, The University of St Andrews, Scotland.\\
	 \vspace{3mm}
	$^\dagger$Department of Pure Mathematics, The University of Waterloo, Canada.
	\end{center}
\begin{abstract}
We consider the Assouad spectrum, introduced by Fraser and Yu, along with a natural
variant that we call the `upper Assouad spectrum'. These spectra are designed  to interpolate
between the upper box-counting and Assouad dimensions. It is known that the Assouad spectrum approaches the upper box-counting dimension at the left hand side of its domain, but does not necessarily approach the Assouad dimension on the right.  Here we show that it necessarily approaches the
\emph{quasi-Assouad dimension} at the right hand side of its domain.
We further show that the upper Assouad spectrum can be expressed in terms of the Assouad
spectrum, thus motivating the definition used by Fraser-Yu.   

We also provide a large family of examples demonstrating new phenomena relating to the form of the Assouad spectrum. 
For example, we prove that it can be strictly concave, exhibit phase transitions of any order, and need not be piecewise differentiable.   \\


\end{abstract}

\section{Assouad type dimensions and spectra}
The Assouad dimension of a metric space is a highly localised measure of its `thickness'. Due to
this it is an important tool when studying bi-Lipschitz embeddings of metric spaces. While the
Assouad dimension captures the worst local covering of a space, its box-counting dimension is a
more `averaged' measure of scaling complexity.
Fraser and Yu introduced the Assouad spectrum as a tool to interpolate between the upper box-counting and
Assouad dimensions, see \cite{Spectraa, Spectrab}.  The Assouad spectrum necessarily approaches the
upper box-counting dimension at the left hand side of its domain, but it was shown in \cite{Spectraa} that
it need not approach the Assouad dimension at the right hand side.  
Similar to the Assouad dimension, the quasi-Assouad dimension is also an upper bound to the Assouad spectrum.
It differs from the Assouad dimension by ignoring `sub-exponential effects'. 
While in most natural settings the quasi-Assouad and Assouad dimensions coincide, the quasi-Assouad
dimension can be strictly smaller.
In those examples where the Assouad spectrum reaches the Assouad dimension, the quasi-Assouad and
Assouad dimensions coincide and it is natural to ask whether the Assouad spectrum always approaches the
quasi-Assouad dimension.
In this article we prove that this is indeed the case and we also exhibit a new family of examples of
possible spectra.

 For a bounded set $E \subseteq \mathbb{R}^d$ and a scale $r>0$ we let $N(E,r)$ be the minimum
 number of sets of diameter $r$ required to cover $E$.  The Assouad dimension of a set $F
 \subseteq \mathbb{R}^d$  is defined by
\[
\dim_\text{A} F \ = \  \inf \Bigg\{ \alpha \  : \ ( \exists  \, C>0) \, (\forall \, 0<r<R) \ \sup_{x \in F} \, N\big( B(x,R) \cap F, r \big) \ \leq \ C \bigg(\frac{R}{r}\bigg)^\alpha  \Bigg\}
\]
and the quasi-Assouad dimension, introduced much more recently by L\"u and Xi \cite{LuXi},  is defined by

\newpage

\begin{multline}\nonumber
\dim_{\mathrm{qA}} F \ = \  \lim_{\theta \to 1} \Bigg(\inf \bigg\{ \alpha \  : \   (\exists C>0) \,
  (\forall 0<r \leq R^{1/\theta} \leq R< 1) \,  (\forall x \in F) \\  N \big( B(x,R) \cap F ,r \big)
  \ \leq \ C \left(\frac{R}{r}\right)^\alpha \bigg\}\Bigg).
\end{multline}

We can see from its definition that the quasi-Assouad dimension leaves an `exponential gap' between $r$ and $R$.
 This gap can be exploited in some
stochastic settings to show that the quasi-Assouad dimension behaves more like the upper
box-counting dimension than the Assouad dimension. Interesting examples of such behaviour are
Mandelbrot and fractal percolation (on self-similar sets): the Assouad dimension is almost surely the  dimension of the percolated
set (which is as big as possible), see  \cite{Fraser14a,Troscheit18}, whereas the quasi-Assouad dimension almost surely coincides with the upper box-counting dimension (which is almost surely strictly smaller than the ambient dimension), see \cite{Spectrab,Troscheit18}. 

The Assouad spectrum, introduced by Fraser and Yu \cite{Spectraa},  is the function defined by
\[
\theta  \ \mapsto \ \dim_{\mathrm{A}}^\theta F \ = \  \inf \bigg\{ \alpha \  : \   (\exists C>0) \, (\forall 0<R<1) \,  (\forall x \in F) \,   N \big( B(x,R) \cap F ,R^{1/\theta} \big) \ \leq \ C \left(\frac{R}{R^{1/\theta}}\right)^\alpha \bigg\}
\]
where  $\theta$ varies over  $(0,1)$. Here, the parameter $\theta$ fixes the relationship
between $r=R^{1/\theta}$ and $R$, but it is equally natural to consider the  `upper spectrum', which 
fixes $R^{1/\theta}$ as an upper bound to $r$ only, defined by
\[
\theta  \ \mapsto \ \overline{\dim}_{\mathrm{A}}^\theta F \ = \  \inf \bigg\{ \alpha \  : \
  (\exists C>0) \, (\forall 0<r \leq R^{1/\theta} \leq R< 1) \,  (\forall x \in F) \,  N \big(
  B(x,R) \cap F ,r \big) \ \leq \ C \left(\frac{R}{r}\right)^\alpha \bigg\},
\]
where, again, $\theta$ varies over  $(0,1)$. We remark that $\overline{\dim}_{\mathrm{A}}^\theta F$ corresponds
to $h_F(\delta)$ in \cite{LuXi}, where $\delta=1/\theta-1$.

We write $\overline{\dim}_\text{B} F$ for the upper box-counting dimension but refer the reader to
\cite{Falconer} for the definition. When we discuss the upper box-counting dimension we are implicitly
referring to bounded sets only, since the definition does not readily apply to unbounded sets. For $F \subseteq \mathbb{R}^d$ and any $\theta \in (0,1)$, we have
\begin{equation}\label{eq:dimensionInequalities}
\overline{\dim}_\text{B} F \leq  \dim_{\mathrm{A}}^\theta F \leq  \overline{\dim}_{\mathrm{A}}^\theta F \leq \dim_{\mathrm{qA}} F  \leq \dim_{\mathrm{A}} F  
\end{equation}
and by definition $ \overline{\dim}_{\mathrm{A}}^\theta F \to  \dim_{\mathrm{qA}} F $ as $\theta \to 1$.  It was also shown in \cite{Spectraa} that $\dim_{\mathrm{A}}^\theta F$ is a continuous function of $\theta$ and  satisfies
\begin{equation} \label{spectrumbound}
 \dim_{\mathrm{A}}^\theta F  \leq \frac{\overline{\dim}_\text{B} F}{1-\theta}
\end{equation}
and therefore  $ \dim_{\mathrm{A}}^\theta F \to  \overline{\dim}_\textup{B} F $ as $\theta \to 0$.  Also, note that by definition $\overline{\dim}_{\mathrm{A}}^\theta F$ is non-decreasing in $\theta$, but it was shown in  \cite[Section 8]{Spectraa} that $\dim_{\mathrm{A}}^\theta F$ is not necessarily non-decreasing: in particular,  $\dim_{\mathrm{A}}^\theta F$ and   $\overline{\dim}_{\mathrm{A}}^\theta F$ do not necessarily coincide.  

For more background on the Assouad dimension, see \cite{Luukkainen, Robinson, Fraser14},  for the
quasi-Assouad dimension, see \cite{LuXi, Hare}, and for the upper box-counting dimension, see \cite{Falconer}.

\section{Results}

Our main technical theorem, which we prove in Section \ref{mainproof}, is the following.

\begin{thm} \label{main}
Let $F \subseteq \mathbb{R}^d$. Then, for all $\theta \in (0,1)$, 
\[
\overline{\dim}_{\mathrm{A}}^\theta F = \sup_{0< \theta' \leq  \theta} \dim_{\mathrm{A}}^{\theta'} F .
\]
\end{thm}

This result shows that all of the information contained in the upper spectrum is also contained in
the Assouad spectrum.  This has the benefit of focusing future study on the Assouad spectrum rather
than the  upper spectrum which could have  \emph{a priori} contained new information in its own
right.  Moreover, as a corollary we obtain the interpolation result which motivated the introduction
of the Assouad spectrum in the first place, albeit with Assouad dimension replaced by quasi-Assouad
dimension.

\begin{cor} \label{maincor}
Let $F \subseteq \mathbb{R}^d$. Then $ \dim_{\mathrm{A}}^\theta F \to
\dim_{\mathrm{qA}}F$ as $\theta \to 1$.
\end{cor}

Theorem \ref{main} only directly implies that $\limsup_{\theta \to 1} \dim_{\mathrm{A}}^\theta F =
\dim_{\mathrm{qA}}F$, but the fact that $\dim_{\mathrm{A}}^\theta F$ has a limit as $\theta \to 1$
follows from estimates in \cite{Spectraa}.  We give the details in Section \ref{corproof}. Combining \eqref{spectrumbound} and Corollary \ref{maincor} we immediately obtain the following
result.

\begin{cor} \label{maincor2}
Let $F \subseteq \mathbb{R}^d$.  Then  $ \overline{\dim}_{\mathrm{B}} F =0$ if and only if
$\dim_{\mathrm{qA}}F=0$.
\end{cor}

The `only if' part of this result is  surprising since one generally has control in the opposite
direction, that is $ \overline{\dim}_{\mathrm{B}} F  \leq \dim_{\mathrm{qA}}F$, and, moreover, the
\emph{Assouad} dimension can take on any value in $[0,d]$ even in cases where the box-counting dimension is
0.  We are not aware of such a `null-equivalence' result holding for any other pair of dimensions.

There are several other consequences of Theorem \ref{main} regarding the upper spectrum which can be
derived from analogous properties of the Assouad spectrum.  For example, the upper spectrum is
immediately seen to be continuous and to approach the upper box-counting dimension as $\theta \to
0$.  It also follows immediately that the two spectra coincide on any interval where the upper
spectrum is strictly increasing.  An example was constructed in  \cite[Section 8]{Spectraa}
demonstrating that the Assouad spectrum can be strictly decreasing (and thus distinct from the upper
spectrum) on infinitely many disjoint intervals accumulating at $\theta = 0$.  However, we believe
this behaviour is not possible at $\theta=1$ and make the following conjecture.

\begin{conj}
For any  set $F \subseteq \mathbb{R}^d$, there exists $\theta_0 \in (0,1)$ such that $\overline{\dim}_{\mathrm{A}}^\theta F =  \dim_{\mathrm{A}}^{\theta} F$ for all $\theta \in [\theta_0,1)$.
\end{conj}

Finally, we present a new family of examples concerning the Assouad spectrum.  In \cite{Spectraa}
some consideration was given to the possible forms the spectrum can take, however, many questions
remain open.  In particular, in all examples so far the spectrum has been piecewise convex,
piecewise analytic, and the only examples of phase transitions have been  points of
non-differentiability.  It follows from \cite[Corollary 3.7]{Spectraa} that the spectrum is
Lipschitz on every closed interval strictly contained in $(0,1)$ and therefore it is differentiable
almost everywhere by Rademacher's Theorem. However, in all examples so far the points of
non-differentiability have been a finite set, or a discrete set accumulating only at 0.  Finally, in
all previous examples the spectrum has been constant in a neighbourhood of $1$.  Here we demonstrate
that much richer behaviour is possible.

\begin{thm} \label{concave}
Let $f:[0,1] \to [0,1]$ be continuous, concave, non-decreasing and satisfy $f(0)>0$ and $f(\theta) \leq  f(0)(\theta+1)$ for all $\theta \in [0,1]$.  Then there exists a compact set $F \subseteq [0,1]$ such that
\[
\dim_{\mathrm{A}}^{\theta} F = f(\theta)
\]
for all $\theta \in (0,1)$.
\end{thm}

The proof of Theorem \ref{concave} will be given in Section \ref{concaveproof}.  The proof gives a recipe for constructing further examples and we have not attempted to optimise its utility for sake of clarity.  The basic strategy is to establish countable stability of the spectrum in a very specific situation and then build the desired function from known examples.   Note that the spectrum is \emph{not} generally countably stable. Theorem \ref{concave} demonstrates that the following list of phenomena are possible, all of which have not been seen before:
\begin{enumerate}
\item The points of non-differentiability can be dense in $(0,1)$.
\item  Phase transitions of all orders are possible, that is points at which the spectrum is $C^k$ but not $C^{k+1}$.
\item  The spectrum need not be constant in a neighbourhood of 1.
\item The spectrum is not necessarily piecewise analytic, or even piecewise differentiable.  This answers \cite[Question 9.1]{Spectraa} in the negative.
\item The spectrum can be strictly concave.
\item The spectrum can be simultaneously strictly increasing and analytic on the whole interval $(0,1)$.
\end{enumerate}

Note that for all examples provided by Theorem \ref{concave}, the upper spectrum and Assouad spectrum coincide since the Assouad spectrum is non-decreasing.

\section{Proofs}

\subsection{A tale of two spectra: proof of Theorem \ref{main}} \label{mainproof}

Let $\theta \in (0,1)$,  suppose $s = \overline{\dim}_\textup{A}^\theta F>0$, and  let $0<\varepsilon< s$. Note that if $s=0$, then the result is trivial.  By definition we can find sequences $x_i, r_i, R_i$ $(i \in \mathbb{N})$ such that $x_i \in F$,  $0<r_i \leq R_i^{1/\theta} < R_i < 1$, $R_i \to 0$ and
\begin{equation} \label{keyest}
N\left( B(x_i,R_i) \cap F, r_i \right) \geq \left( \frac{R_i}{r_i} \right)^{s-\varepsilon}.
\end{equation}
For each $i$, let $\theta_i$ be defined by $r_i = R_i^{1/\theta_i}$, noting that $0 < \theta_i \leq \theta$.  Using compactness of $[0,\theta]$ to extract a convergent subsequence, we may assume that $\theta_i \to \theta' \in [0,\theta]$ and by taking a further subsequence if necessary we can assume that $|\theta_i - \theta' | < \delta$ for all $i$ where $\delta >0$ can be chosen arbitrarily. We may also assume that the sequence $\theta_i$ is either non-increasing or strictly increasing.   Assume for now that $\theta'>0$.  We will deal with the $\theta'=0$ case separately at the end.

If the sequence $\theta_i$ is non-increasing, then $\theta' \leq \theta_i$, and therefore $R_i^{1/\theta_i} \geq R_i^{1/\theta'}$, for all $i$.  It follows that
\begin{eqnarray*}
N\left( B(x_i,R_i) \cap F, R_i^{1/\theta'} \right)  \geq N\left( B(x_i,R_i) \cap F, R_i^{1/\theta_i} \right)
 &\geq & \left( \frac{R_i}{R_i^{1/\theta_i}} \right)^{s-\varepsilon} \qquad \text{by \eqref{keyest}} \\ \\ &=&   \left( \frac{R_i}{R_i^{1/\theta'}} \right)^{\left(\frac{1-1/\theta_i}{1-1/\theta'}\right)(s-\varepsilon)} \\ \\
 &\geq&\left( \frac{R_i}{R_i^{1/\theta'}} \right)^{\frac{\theta'(1-\theta'-\delta)}{(\theta'+\delta)(1-\theta')}(s-\varepsilon)}
\end{eqnarray*}
This yields $\dim_\textup{A}^{\theta'}F \geq \frac{\theta'(1-\theta'-\delta)}{(\theta'+\delta)(1-\theta')}(s-\varepsilon)$ and, since $\delta >0$ can be chosen arbitrarily small (after fixing $\theta'$), we obtain $\dim_\textup{A}^{\theta'}F \geq s-\varepsilon$.

On the other hand, if $\theta_i$ is strictly increasing, then $\theta' > \theta_i$ for all $i$.
Taking another subsequence if necessary we can also assume that $\theta_i > \theta'/2$ for all $i$.
Covering by $R_i^{1/\theta'}$-balls and then covering each $R_i^{1/\theta'}$-ball by
$R_i^{1/\theta_i}$-balls we obtain
\begin{eqnarray*}
N\left( B(x_i,R_i) \cap F, R_i^{1/\theta_i} \right) &\leq& N\left( B(x_i,R_i) \cap F, R_i^{1/\theta'} \right)\left( \sup_{z \in \mathbb{R}^d} N\left(B(z,R_i^{1/\theta'}), R_i^{1/\theta_i} \right) \right) \\ \\
& \leq& N\left( B(x_i,R_i) \cap F, R_i^{1/\theta'} \right) c(d) \left( \frac{R_i^{1/\theta'}}{R_i^{1/\theta_i}} \right)^{d}
\end{eqnarray*}
where $c(d)\geq 1$ is a constant depending only on the ambient spatial dimension  $d$.  Therefore
\begin{eqnarray*}
N\left( B(x_i,R_i) \cap F, R_i^{1/\theta'} \right)  &\geq& c(d)^{-1} N\left( B(x_i,R_i) \cap F, R_i^{1/\theta_i} \right)  R_i^{(1/\theta_i-1/\theta')d} \\ \\
& \geq& c(d)^{-1} \left( \frac{R_i}{R_i^{1/\theta_i}} \right)^{s-\varepsilon}
R_i^{(1/\theta_i-1/\theta')d}  \qquad \text{by \eqref{keyest}} \\ \\
& =& c(d)^{-1} \left( \frac{R_i}{R_i^{1/\theta'}} \right)^{\left(\frac{1-1/\theta_i}{1-1/\theta'}\right)(s-\varepsilon)+\left(\frac{1/\theta_i-1/\theta'}{1-1/\theta'}\right)d} \\ \\
& \geq & c(d)^{-1} \left( \frac{R_i}{R_i^{1/\theta'}} \right)^{s-\varepsilon-\frac{\delta d}{(1-\theta')\theta'/2}}.
\end{eqnarray*}
It follows that  $\dim_\textup{A}^{\theta'}F \geq s-\varepsilon-\frac{\delta d}{(1-\theta')\theta'/2}$ and since $\delta >0$ can be chosen arbitrarily small (after fixing $\theta'$) we obtain $\dim_\textup{A}^{\theta'}F \geq s-\varepsilon$ as before.  Since $\varepsilon>0$ was arbitrary it follows that
\[
\sup_{0< \theta' \leq \theta} \dim_\textup{A}^{\theta'}F \geq s
\]
completing the proof, noting that the other direction is trivial.

All that remains is to consider the case where $\theta'=0$.  Interestingly, this case is very straightforward if $F$ is bounded, but not otherwise.  Indeed, for bounded $F$,
\[
 N\left( F, R_i^{1/\theta_i} \right) \geq N\left( B(x_i,R_i) \cap F, R_i^{1/\theta_i} \right) \geq  \left( \frac{R_i}{R_i^{1/\theta_i}} \right)^{s-\varepsilon} \geq   \left( \frac{1}{R_i^{1/\theta_i}} \right)^{(s-\varepsilon)(1-\delta)}
\]
by \eqref{keyest}.  Then, by \eqref{eq:dimensionInequalities}, $\dim_\textup{A}^\theta F \geq
\overline{\dim}_\textup{B}F \geq (s-\varepsilon)(1-\delta)$.  Since $\delta>0$ and  $\varepsilon>0$
can be chosen arbitrarily small, this yields the desired result.  However, if $F$ is unbounded, then we cannot easily go via box-counting dimension and more work is needed.  Let $\phi \in (0, \theta)$ be chosen such that
\[
\frac{\log \theta}{\log \phi} \notin \mathbb{Q}
\]
and suppose for a contradiction that $\max\{ \dim_\textup{A}^{\phi}F, \dim_\textup{A}^{\theta}F\} < s-2 \varepsilon$.  It follows that for all $R$ small enough and all $x \in F$  we have
\begin{equation} \label{cover11}
 N\left( B(x,R), R^{1/\phi} \right) \leq  \left( \frac{R}{R^{1/\phi}} \right)^{ s-2 \varepsilon}
\end{equation}
and
\begin{equation} \label{cover22}
 N\left( B(x,R), R^{1/\theta} \right) \leq  \left( \frac{R}{R^{1/\theta}} \right)^{ s-2 \varepsilon}.
\end{equation}
Note that we can get rid of any constants here since we are only considering the spectrum at two points, and therefore two instances of having to take small enough $R$.  Consider the additive monoid generated by $\{\log \phi, \log \theta\}$, that is, the set
\[
\{ m \log \phi + n \log \theta \ : \ m,n \geq 0, \, m,n \in \mathbb{Z}\} \subset (-\infty, 0].
\]
By our irrationality assumption on $\phi$ and $\theta$, it follows that for all $\eta>0$ if $i$  is
large enough there always  exists $m,n \geq 0$ such that
\[
0 \leq \log (\phi^m \theta^n) - \log \theta_i \leq \eta.
\]
In particular, this implies that  for sufficiently large $i$ we can choose $m,n \geq 0$ such that
\begin{equation} \label{goodchoice}
0 \leq 1/\theta_i - 1/(\phi^{m}\theta^n) \leq \varepsilon/(2d\theta_i).
\end{equation}
Fix a large $i$ and $m,n$ corresponding to $i$ as in \eqref{goodchoice} above.  We can now build an
efficient cover of $B(x_i,R_i)$ by $R_i^{1/\theta_i}$-balls.  We begin by covering $B(x_i,R_i)$ with
$R_i^{1/\phi}$-balls and then each of these $R_i^{1/\phi}$-balls by $R_i^{1/\phi^2}$-balls and
continue in this way until we have covered $R_i^{1/\phi^{m-1}}$-balls by $R^{1/\phi^m}$-balls.  We
then switch to a `$\theta$-regime', covering each $R_i^{1/\phi^m}$-ball with
$R_i^{1/(\phi^m\theta)}$-balls.  Each of these balls is covered by $R_i^{1/(\phi^m\theta^2)}$-balls
until we reach a covering by $R_i^{1/(\phi^m\theta^n)}$-balls.  Using \eqref{cover11}-\eqref{cover22} and a standard telescoping argument we therefore get
\begin{eqnarray*}
 N\left( B(x_i,R_i), R_i^{1/(\phi^m\theta^n)} \right) &\leq&  \left( \frac{R_i}{R_i^{1/\phi}} \right)^{ s-2 \varepsilon}\left( \frac{R_i^{1/\phi}}{R_i^{1/\phi^2}} \right)^{ s-2 \varepsilon}  \cdots  \left( \frac{R_i^{1/\phi^{m-1}}}{R_i^{1/\phi^{m}}} \right)^{ s-2 \varepsilon} \\ \\
&\,& \qquad \times  \left( \frac{R_i^{1/\phi^{m}}}{R_i^{1/(\phi^{m}\theta)}} \right)^{ s-2
\varepsilon}\left( \frac{R_i^{1/(\phi^{m}\theta)}}{R_i^{1/(\phi^{m}\theta^2)}} \right)^{ s-2
\varepsilon}  \cdots  \left( \frac{R_i^{1/(\phi^{m}\theta^{n-1})}}{R_i^{1/(\phi^{m}\theta^n)}} \right)^{ s-2 \varepsilon} \\ \\
&=&  \left( \frac{R_i}{R_i^{1/(\phi^{m}\theta^n)}} \right)^{ s-2 \varepsilon}.
\end{eqnarray*}
Finally, to obtain a cover by $R_i^{1/\theta_i}$-balls we cover each $R_i^{1/(\phi^{m}\theta^n)}$-ball by at most
\[
c(d) \left( \frac{R_i^{1/(\phi^{m}\theta^n)}}{R_i^{1/\theta_i}} \right)^{ d}
\]
many $R_i^{1/\theta_i}$-balls, where $c(d)\geq 1$ is, as above,  a constant depending only on the ambient spatial
dimension $d$. Combining this with \eqref{keyest} we get that for all large enough $i$ (and thus
small enough $R_i$) we must have
\[
 \left( \frac{R_i}{R_i^{1/\theta_i}} \right)^{s-\varepsilon} \leq N\left( B(x_i,R_i) \cap F, r_i
 \right) \leq   c(d) \left( \frac{R_i}{R_i^{1/(\phi^{m}\theta^n)}} \right)^{ s-2 \varepsilon}  \left(
 \frac{R_i^{1/(\phi^{m}\theta^n)}}{R_i^{1/\theta_i}} \right)^{ d}.
\]
Using \eqref{goodchoice} we therefore have
\[
(1/\theta_i -1) (s-\varepsilon) \leq \left(1/(\phi^{m}\theta^n\right) - 1)(s-2 \varepsilon) + (1/\theta_i -
1/(\phi^{m}\theta^n)) d \leq  (1/\theta_i - 1)(s-2 \varepsilon) + \varepsilon/(2\theta_i) 
\]
which, since $\theta_i \to 0$, yields $s-\varepsilon \leq s- 3 \varepsilon/2$, a contradiction.  It
follows that the Assouad spectrum is at least $s-2 \varepsilon$ at either $\theta$ or $\phi$, which
upon letting $\varepsilon \to 0$ proves the result.\qed

\subsection{Interpolation in the limit: proof of Corollary \ref{maincor}} \label{corproof}

As already stated, Theorem \ref{main}  directly implies that $\limsup_{\theta \to 1} \dim_{\mathrm{A}}^\theta F = \dim_{\mathrm{qA}}F$.  Therefore all that remains is to prove that $\lim_{\theta \to 1} \dim_{\mathrm{A}}^\theta F$ exists.  This follows from results in \cite[Section 3]{Spectraa}, although it was not explicitly stated.  In particular, we have the following lemma, which appears as part of \cite[Remark 3.9]{Spectraa}.

\begin{lma} \label{nthroot}
For nonempty $F \subseteq \mathbb{R}^d$, $\theta \in (0,1)$ and $n \in \mathbb{N}$, we have
\[
\dim_{\mathrm{A}}^\theta F \leq \dim_{\mathrm{A}}^{\sqrt[n]{\theta}} F.
\]
\end{lma}

Let $t = \limsup_{\theta \to 1} \dim_{\mathrm{A}}^\theta F = \dim_{\mathrm{qA}}F$ and $\varepsilon>0$.  Since $\dim_{\mathrm{A}}^\theta F $ is a continuous function of $\theta$ we can find $0<a<b<1$ such that for all $\theta \in [a,b]$ we have $\dim_{\mathrm{A}}^\theta F > t-\varepsilon$.  It follows from Lemma \ref{nthroot} that for all
\[
\theta \in \bigcup_{n \in \mathbb{N}} \left[\sqrt[n]{a}, \sqrt[n]{b} \right] =: X
\]
we also have  $\dim_{\mathrm{A}}^\theta F > t-\varepsilon$.  However, it is easily seen that $X$ contains an interval $(x,1)$ for some $x \in (0,1)$.  Indeed, the intervals $[\sqrt[n]{a}, \sqrt[n]{b} ]$  and $[\sqrt[n+1]{a}, \sqrt[n+1]{b} ]$ intersect each other when $n$ is sufficiently large.  In fact, one can choose $x=\sqrt[n]{a}$ where $n$ is chosen large enough to ensure that
\[
\frac{n}{n+1} \geq \frac{\log b}{\log a}.
\]
It follows that $\liminf_{\theta \to 1} \dim_{\mathrm{A}}^\theta F = \limsup_{\theta \to 1}
\dim_{\mathrm{A}}^\theta F = t$, as required.\qed

\subsection{New examples: proof of Theorem \ref{concave}} \label{concaveproof}

The Moran constructions considered in \cite{Spectrab} provide us with a simple but useful family of examples.  In particular, we have the following result by applying \cite[Corollary 6.2]{Spectrab} to the examples considered towards the end of \cite[Section 6.2]{Spectrab}.

\begin{lma}[Section 6.2, \cite{Spectrab}]
For any $0<s<t \leq 1$, there exists a compact set $F \subseteq [0,1]$ such that
\[
\dim_{\mathrm{A}}^\theta F = \min \left\{\frac{s}{1-\theta}, \, t \right\}.
\]
\end{lma}

Note that these examples attain the general upper bound \eqref{spectrumbound} until the
quasi-Assouad dimension is reached. Such sets $F$ are constructed in \cite{Spectrab}  as homogeneous
(dyadic) Moran constructions where one has complete control over the number of dyadic intervals
present inside a higher level dyadic interval.  Thus one also has complete control, up to a uniform
constant, on the covering numbers $N(B(x,R) \cap F , R^{1/\theta})$.  Therefore, either following
the proof of \cite[Corollary 6.2]{Spectrab} or simply `pruning' the  sets $F$ as necessary, one can
`upgrade' the above lemma as follows.

\begin{lma} \label{exampless}
For any $0<s<t \leq 1$, there exists a compact set $F \subseteq [0,1]$ such that
\[
\dim_{\mathrm{A}}^\theta F = \min \left\{\frac{s}{1-\theta}, \, t \right\} =: u(\theta)
\]
and, moreover, for all $x \in F$ and $R \in (0,1)$ we have
\[
N(B(x,R) \cap F , R^{1/\theta}) \leq 10 \left( \frac{R}{R^{1/\theta}} \right)^{u(\theta)}.
\]
\end{lma}

Let $\{q_i\}_{i \geq 1}$ be an enumeration of the rationals in $(0,1)$ and for each $i$, let $F_i \subseteq [0,1]$ be the set provided by Lemma \ref{exampless} where $t=t_i = f(q_i)$ and $s=s_i=f(q_i) (1-q_i)$.   In particular the phase transition in $\dim_{\mathrm{A}}^\theta F_i  =: u_i(\theta)$ occurs with coordinates $(q_i, f(q_i))$.  Also, note that by assumption
\[
f(\theta) \leq f(0)(\theta+1) \leq \frac{f(0)}{1-\theta}
\]
and therefore $s_i \leq f(0)$ for all $i$. Since $f$ is concave and non-decreasing and $u_i$ is
convex on $[0,q_i]$ it follows that
$u_i(\theta) \leq f(\theta)$ for all $\theta \in (0,1)$ and that $u_i(q_i) = f(q_i)$ for all $i \geq
1$.  Therefore, since $f$ is continuous, we can conclude that
\[
\sup_{i \geq 1} u_i(\theta) = f(\theta)
\]
for all $\theta \in (0,1)$.

We can now construct the set $F$ required to prove the theorem.  Let
\[
F = \{0\} \cup \bigcup_{i \geq 1} \widehat{F}_i
\]
where $\widehat{F}_i=2^{-2^i} F_i + 2^{-2^i} = \{ 2^{-2^i} x +2^{-2^i}: x \in F_i\} \subseteq
[0,1]$.  Let $\theta \in (0,1)$, $x \in F$ and $R \in (0,1)$.  Let $j = \min\{ i \geq 1 :
  \widehat{F}_i \cap
B(x,R) \neq \emptyset\}$.  If $\widehat{F}_{j'} \cap B(x,R) \neq \emptyset$ for some $j' > j$, then $R \geq
2^{-2^j}/4$ and also $R^{1/\theta} \geq (2^{-2^j/\theta})/4^{1/\theta}$.  Therefore, there is a
constant $k \in \mathbb{N}$ depending only on $\theta$ such that
\[
N\left(\{0\} \cup \bigcup_{i \geq j+k} \widehat{F}_i , \ R^{1/\theta}\right) \leq 1.
\]
Therefore by Lemma \ref{exampless} we conclude that
\[
N(B(x,R) \cap F , R^{1/\theta}) \leq 1+ \sum_{i=j}^{j+k-1} N(B(x,R) \cap F_i , R^{1/\theta}) \leq 
1+ 10 k  \left( \frac{R}{R^{1/\theta}} \right)^{\sup_{i \geq 1}u_i(\theta)}.
\]
This proves that $\dim_{\mathrm{A}}^\theta F  \leq \sup_{i \geq 1} u_i(\theta) = f(\theta)$.  The
reverse inequality is immediate by monotonicity and therefore the theorem is proved.\qed

\section*{Acknowledgements}

This work began whilst JMF visited the University of Waterloo in March 2018. He is grateful for the
financial support, hospitality, and inspiring research atmosphere.  JMF was financially supported by a
\emph{Leverhulme Trust Research Fellowship} (RF-2016-500) and  an \emph{EPSRC Standard Grant}
(EP/R015104/1).  
KEH was supported by \emph{NSERC Grant} 2016-03719. 
KGH was supported by \emph{NSERC Grant} 2014-03154.
ST was supported by \emph{NSERC Grants} 2016-03719 and 2014-03154, and the University of Waterloo.
HY was financially supported by the University of St Andrews.\\

\vskip.5cm
\noindent
---------------------

{\small
Jonathan M. Fraser, E-mail: \texttt{jmf32@st-andrews.ac.uk}
\vskip.2em
Kathryn E. Hare, E-mail:  \texttt{kehare@uwaterloo.ca}
\vskip.2em
Kevin G. Hare, E-mail:  \texttt{kghare@uwaterloo.ca}
\vskip.2em
Sascha Troscheit, E-mail:  \texttt{stroscheit@uwaterloo.ca}
\vskip.2em
Han Yu, E-mail: \texttt{hy25@st-andrews.ac.uk}
}
\end{document}